\documentclass[12pt]{article} 
\usepackage[french]{babel}
\usepackage[latin1]{inputenc}
\usepackage{amsfonts}
\usepackage{latexsym}

\newtheorem{thm}{Th\'eor\`eme}[section]
\newtheorem{prop}[thm]{Proposition}
\newtheorem{lem}[thm]{Lemme}

\newtheorem{defi}[thm]{D\'efinition}
\newtheorem{rmq}[thm]{Remarque}

\newenvironment{preuve}{\par\indent\textit{Preuve.\/}}{\raisebox{-.2ex}{$\Box$}}

\newcommand{\carre}{\raisebox{-.2ex}{$\Box$}}
\newcommand{\C}{\ensuremath{\mathbb C}}
\newcommand{\R}{\ensuremath{\mathbb R}}

\newcommand{\N}{\ensuremath{\mathbb N}}

\begin{document}

\title{Forme normale formelle d'une perturbation à séparatrices fixées d'un champ
hamiltonien quasi-\-ho\-mo\-gène }
\author{\sc E. Paul}
\date{\today}
\maketitle

%\begin{center}
%\textit{Version préliminaire FNF.5}
%\end{center}

{\bf Résumé.} {\footnotesize On classifie à changement de variable formel près
les champs de vecteurs --ainsi que les feuilletages qu'ils définissent-- qui
sont des perturbations d'un champ hamiltonien quasi-homogène $X_0$ par des
termes de degré de quasi-homogénéité supérieur. Le degré $\delta_0$ du champ
$X_0$ est quelconque, mais on demande que le champ perturbé laisse encore
invariant les séparatrices de $X_0$. Les formes normales présentées ici
généralisent celles obtenues lorsque $X_0$ possède une partie linéaire
semi-simple ou nilpotente. Nous les interprétons géométriquement en terme de
cohomologie relative au champ initial.}

\bigskip
\bigskip

\noindent{\Large \textbf{Introduction}}

\bigskip

On s'intéresse à la classification formelle de germes de champs de vecteurs
analytiques à l'origine de $\C^n$ qui sont des perturbations d'un champ de
vecteurs quasi-\-ho\-mo\-gène $X=X_0+\cdots$, les termes perturbateurs ayant
des degrés de quasi-homogénéité supérieur à la partie initiale $X_0$. Deux
relations d'équivalences nous intéressent ici:

\bigskip

i) \textit{Classification des champs:} deux champs $X$ et $Y$ sont formellement
(resp. analytiquement) conjugués ($X \sim Y$) s'il existe une transformation formelle
(resp. analytique) $\Phi $ telle que
$Y=\Phi ^*X.$

\bigskip 

ii) \textit{Classification des feuilletages:} deux champs $X$ et $Y$ sont
formellement (resp. analytiquement) orbitalement équivalents ($X
\approx Y$) s'il existe une transformation formelle (resp. analytique) $\Phi $
et une unité formelle (resp. analytique) $u$ telles que
$Y=u\Phi ^*X.$ Dans ce cas, la conjugaison $\Phi $ envoie les orbites du
premier champ sur celles du second, sans nécessairement conjuguer leurs flots.

\bigskip 

Le problème des formes normales, c'est-à-dire de la détermination d'un
repré\-sen\-tant privilégié dans les classes d'équivalences modulo $\sim$ ou
$\approx$, est motivé d'une part par la nécessité de classifier ces champs ou
feuilletages, d'autre part par l'obtention de modèles sur lequel on pourra
répondre à des problèmes ''d'intégrabilité'': existences de solutions ou
intégrales premières dans une classe de transcendance donnée. Pour cette
seconde motivation, l'obtention d'un représentant unique dans uns classe
d'équivalence donnée n'est pas toujours nécessaire. La détermination explicite
de formes normales a été très étudiée lorsque le premier jet du champ est non
nul, ce que nous ne supposerons plus ici.

\bigskip 

Sur $\C^2$ muni de la forme volume standard, à tout champ $X$ on peut associer
la 1-forme duale $\omega _X=i_Xdx\wedge dy.$ Les relations $\sim$ et $\approx$
se définissent de manière analogue pour les 1-formes holomorphes. D'après la formule $\Phi ^*\omega
_X=\mbox{Jac}(\Phi )\omega _{\Phi ^*X}$, on a
$$\omega _X\approx\omega _Y\Leftrightarrow X\approx Y$$ ce qui nous permet de
travailler indifféremment avec les champs ou les formes pour cette
relation. Il n'en est pas de même pour la relation de conjugaison $\sim$: les
champs $\frac{\partial}{\partial x}$ et $(1+x)\frac{\partial}{\partial x}$ sont
conjugués d'après le théorème de redressement des champs réguliers, mais leurs
formes duales ne le sont pas puisque $dy$ est fermée alors que $(1+x)dy$ ne l'est pas.

Nous nous intéressons ici essentiellement à la classification locale des
feuilletages en dimension deux, mais obtiendrons au passage un résultat
partiel de classification des champs de vecteurs lorsque le degré de quasi-homogénéité $\delta_0$
du champ initial est nul: ceci concerne en particulier les champs de type
''Poincaré-Dulac'': $X  =\frac{1}{2pq} (qx\partial_x-py\partial_y)+\cdots$
et les champs de type ''noeud-col'': $X=y\partial_y +\cdots$.

\bigskip

La classification formelle des champs de vecteurs en toute dimension repose
sur un résultat préliminaire qui est bien connu lorsque la partie initiale
$X_0$ est un champ semi-simple \cite{MAR}. Nous le généralisons ici pour une
partie initiale quasi-homogène de degré quelconque: 

\bigskip

{\bf Lemme de pré\-nor\-ma\-li\-sation formelle des champs.} \textit{On se donne un
sup\-plé\-men\-taire $V$ de l'image de l'opérateur ad$(X_0)=[X_0,\cdot]$ dans 
un module de champs formels $M$ gradué par un degré de quasi-homogénéité. 
Le champ
$X=X_0+\cdots$ est formellement conjugué à un champ $X_0+Y$ où $Y$ appartient à
$V$.}

\bigskip

Tout espace de formes prénormales $V$ est donc, d'après ce lemme, isomorphe au quotient
$H_K^1=M/B_K^1$, $B_K^1=\mbox{Im(ad$X_0$)}$, l'indice ${\scriptstyle K}$
faisant ici référence aux complexes de Koszul construits à partir de tels opérateurs ad$(X_0)=[X_0,\cdot]$. 
Dans le cadre de la relation d'équivalence des feuilletages, nous pouvons diminuer la taille du sous-espace $V$:

\bigskip

{\bf Lemme de pré\-nor\-ma\-li\-sation formelle des feuilletages.} \textit{On
  se donne un sup\-plé\-men\-taire
$W$ de $A_K^1=\mbox{Im(ad$X_0$)}+\widehat{\mathcal O}_nX_0$ dans $M$. 
Il existe un champ $Y$ appartenant à $W$, un difféo\-mor\-phisme formel $\Phi $ et une unité formelle $u$ tels que
$$\Phi ^*X=u(X_0+Y).$$}

Ces supplémentaires étant fixés, il n'y a pas unicité de ces formes normales
formelles dans $V$ ou $W$, d'où l'appellation
''pré\-nor\-ma\-li\-sation''. Ceci nécessitera une réduction
supplé\-men\-taire que nous ferons dans une seconde étape. Néanmoins, cette étape de prénormalisation doit être
considérée comme l'étape essentielle du processus de normalisation et non comme
une simple opération préliminaire. C'est en particulier elle qui est
susceptible de créer de la divergence dans les conjuguantes et dans les formes
normales elle-mêmes.

\bigskip

Ces lemmes ramènent donc l'étape de pré\-nor\-ma\-li\-sation des champs (resp. des
feuilletages) à  un choix pertinent pour ce sup\-plé\-men\-taire $V$ (resp. $W$), c'est-à-dire
à un choix d'éléments du module $M$ dont la classe engendre le quotient
$H_K^1=M/B_K^1$, (resp. $G_K^1=M/A_K^1$). On peut remarquer que ces quotients
sont des modules sur l'anneau $\mathcal I$ des intégrales premières du champ
initial $X_0$. Il est donc naturel d'imposer lors du choix de $V$ ou $W$ de
préserver cette structure. L'isomorphisme que nous construirons entre $V$ et 
$H_K^1$ (ou entre $W$ et $G_K^1$) sera un isomorphisme de $\mathcal
I$-modules.
Le résultat présenté ici se limite au cadre suivant:

\bigskip

1-  Nous nous plaçons en dimension 2. Le champ quasi-homogène initial $X_0$ est alors
tangent à une hypersurface $S$ (séparatrice) d'équation réduite $h_0=0$. 
Nous supposons que le champ perturbé $X=X_0+\cdots$ reste tangent à cette
hypersurface $S$. En d'autres termes, nous considérons une perturbation de
$X_0$ \textit{à séparatrice fixée}. 
Dans certains cas, (par exemple $h_0=xy$, $h_0=y^2-x^3$), cette hypothèse n'en
est pas véritablement une: la séparatrice persiste et, étant rigide, on la
redresse sur la séparatrice initiale par un changement analytique de
coordonnées.  

Fixer les séparatrices nous conduit à travailler dans le module $M=\widehat\chi(\log S)$ des champs formels
\textit{logarithmiques} pour cette hypersurface $S$, c'est-à-dire des champs
$X$ tels que $X(h_0)$ appartient à l'idéal $(h_0)$. En dimension 2, ce
module est toujours un module libre de rang 2. 
Soit $R$ le champ
$p$-radial (pondéré par les poids $p_i$ de quasi\--homo\-gé\-néité) pour
lequel $X_0$ est quasi-homogène. 
Supposons que les champs $X_0$ et $R$ forment une base de 
$\widehat\chi(\log S)$. (Ce sera par exemple le cas sous l'hypothèse (2-)
ci-dessous). Tout champ logarithmique se décompose alors en une composante \textit{''intégrable''} (en $X_0$)
et une composante
\textit{''dissipative''} (en $R$): la première
garde $h$ comme intégrale première, alors que la forme logarithmique duale de
la seconde n'est jamais fermée. En particulier, en comparant les deux lemmes de
pré\-nor\-malisation ci-dessus, on constate immédiatement que si $V$ est un espace de formes
prénormales pour un champ $X$, sa composante dissipative $W$ est un espace de
formes prénormales pour le feuilletage défini par $X$. L'utilisation de cette
base nous permet de démontrer facilement (voir énoncé \ref{PFNF} ) le

\bigskip

{\bf Théorème 1.}
\textit{Soit $X_0$ un champ quasi-homogène de degré $\delta_0$ par rapport à
  un champ quasi-radial $R$, laissant invariant une courbe $S$. On suppose que les champs
$X_0$ et $R$ forment une base des champs logarithmiques pour $S$.
Soit $X=X_0+\cdots$ une perturbation de $X_0$ laissant fixe $S$.
\begin{enumerate}
\item 
Le module des formes prénormales de feuilletages définis par
$X$ est 
 $$G_K^1=\mbox{Coker}(X_0)\cdot R.$$  
\item Si le degré de quasi-homogénéité $\delta_ 0$ du champ initial est nul,
 le module des formes prénormales du champ $X$ est
$$H_K^1=\mbox{Coker}(X_0)\cdot X_0+\mbox{Coker}(X_0)\cdot R.$$
\end{enumerate}}

\bigskip

Le problème de la détermination d'un espace de formes prénormales est donc
maintenant ramené à la détermination du conoyau de $X_0$.

\bigskip

2-  Nous supposons de plus que le champ initial 
$X_0$ est hamiltonien, dual pour la forme $dx\wedge dy$ d'une forme
$\delta^{-1}dh$, où l'intégrale première $h$ est quasi-homogène de degré
$\delta$. On peut alors vérifier que les champs $X_0$ et $R$ forment une base
de $\widehat\chi(\log S)$. 
Supposons de plus que $h$
soit à singularité isolée.
On désigne par $J_0$
l'idéal des composantes de $X_0$ dans la base usuelle (qui est aussi l'idéal jacobien de $h$), et par $a_1=1,\cdots
a_{\mu}$, une base de monômes engendrant l'espace vectoriel ${\cal O}^2/J_0$.
Nous démontrons (théorème \ref{conoyau})
que le conoyau de la dérivation hamiltonienne $X_0$ est un $\C[[h]]$-module libre de
rang $\mu$ engendré par les champs $a_iR$. Ce fait résulte essentiellement de
l'existence d'une connexion singulière sur cet espace de type Gauss-Manin, et
d'un théorème d'indice de B. Malgrange. De ces deux résultats nous déduisons le

\bigskip

{\bf Théorème 2.}
\textit{Soit $X=X _0 +\cdots $ une perturbation à séparatrices fixées du champ hamiltonien $X _0$. Il
  existe un élément $(d_1,\cdots d_{\mu })$ de $\C[[h]]^{\mu}$, une
  conjugaison formelle $\Phi $ conjuguant orbitalement $X$ au champ formel
$$X _0 +\sum_{i=1}^{\mu }d_i(h){a_iR}.$$ 
De plus, on peut imposer à cette conjugaison d'être fibrée pour le champ $R$, c'est-à-dire d'être l'exponentielle d'un champ colinéaire à $R$.}

\bigskip

Nous prouvons de plus que les coefficients formels $d_i(h)$ des formes prénormales du théorème 2 peuvent s'exprimer par des formules intégrales. Pour
cela, nous considérons le point de vue dual concernant ces formes normales
formelles en les transférant par une ''forme volume logarithmique''
sur le module des formes logarithmiques. Il apparait alors que le module des
formes prénormales de \textit{feuilletages} est isomorphe au module de
cohomologie relative formelle à la forme initiale fermée $\omega_0$
(proposition \ref{calcul.Hrel}). Ceci nous permet, en utilisant une base de formes
horizontales de ce module, d'écrire les coefficients multivalués apparaissant
dans cette base comme des intégrales de la partie dissipative de la
forme normale formelle sur les cycles évanescents de sa partie hamiltonienne.

\bigskip

La non unicité des formes prénormales de feuilletages obtenues au théorème 2 provient d'une
ambiguité dans le choix de la conjuguante fibrée. En effet, on peut faire
opérer sur l'ensemble des formes prénormales d'un feuilletage le groupe des
conjuguantes de la forme $\Phi =\exp b(h)R$ où
$b$ est une série formelle d'une variable. La deuxième étape (réduction
finale) consiste donc à choisir un représentant dans chaque orbite de cette
action sur l'espace des formes prénormales de $X$. Il apparait que, après division par une puissance fractionnaire convenable de $h$, les
coefficients $d_i(h)$ se comportent sous cette action comme des \textit{champs de
vecteurs sur un revêtement fini du disque épointé image de $h$}. Ceci permet de normaliser un des
coefficients  $d_i(h)$ choisi arbitrairement sous forme
$$\frac{ {h}^{m}}{1+\lambda  {h}^{m+n}}.$$
De plus, pour que le résidu $\lambda $ soit non nul, il est nécessaire qu'une
certaine condition de divisibilité soit satisfaite. On obtient ainsi une classification formelle
de ces feuilletages généralisant celles connues par Poincaré et
Dulac, ou celles de F. Loray dans le cas des singularités cuspidales.

L'observation ci-dessus conduisant à cette réduction finale appelle
une interprétation géométrique. Quelle est la signification de ces $\mu$ champs
de vecteurs? On sait d'autre part d'après \cite{JFM}, ou d'après \cite{CM} dans le cas
du cusp, que ce type de
feuilletage avec séparatrices quasi-homogènes est complètement classifié par sa
structure transverse. Quel lien existe-t-il entre ces invariants explicites et
l'holonomie du diviseur exceptionnel considérée dans \cite{CM}?

\bigskip

Remarquons enfin que le processus de réduction finale
ne permet de normaliser sous forme rationnelle qu'une
seule des $\mu $ séries $d_i(h)$ apparaissant dans les formes prénormales. En
conséquence, dès que $\mu $ est plus grand que 1, subsistent des séries
formelles dans ces formes normales. Supposons maintenant que le champ $X$ soit
analytique. La question de la convergence de ces séries se pose alors. On
connait aujourd'hui un exemple (voir \cite{C-D}) où, dans le cas d'une séparatrice de type cusp,
ces séries divergent.
De plus, dans un travail récent, M. Canalis-Durand et R. Schaefke ont établi le caractère 1-sommable des formes normales et
conjuguantes obtenues. Il est donc raisonnable d'espérer une $k$-sommabilité
des formes normales formelles obtenues ici. L'interprétation géométrique
des invariants souhaitée ci-dessus pourrait nous y aider. Nous nous proposons
de développer
ces deux questions dans un travail ultérieur.

\bigskip

Remerciements à Frank Loray --son article \cite{LOR} a été la principale
source d'inspiration de ce travail--, et Reinhart Schaefke,  pour les
discussions que nous avons eu sur ce sujet. 

\pagebreak

\tableofcontents

\bigskip
\bigskip
\bigskip

\noindent {\Large \textbf{Notations}}

\bigskip

On désigne par ${\cal O} _n$ l'anneau des germes de fonctions holomorphes à
l'origine de $\C^n$ et par $\widehat{\cal O} _n$ son complété formel. On
notera $\partial_x$ la dérivation $\frac{\partial}{\partial x}$ et $h_x$ pour
$\frac{\partial h}{\partial x}$.
Dans tout le texte nous noterons:

\begin{itemize}
\item $d_0$: le degré de quasi-homogénéité de l'équation réduite $h_0$ de $S$;
\item $\delta_0$: le degré de quasi-homogénéité du champ initial $X_0$;
\item $\delta$: le degré de quasi-homogénéité de son intégrale première $h$.
\end{itemize}
\pagebreak

\section{Fonctions et champs quasi-\-ho\-mo\-gènes}

Soit $p=(p_1,\cdots p_n)$ une collection
d'entiers positifs tel que les $p_i$ soient premiers entre eux. 
On considère le champ $p$-radial
$$R_{p}=\sum_{i=1}^{n}p_ix_i\frac{\partial}{\partial x_i}.$$
Pour tout monôme $x^m=x_1^{m_1}\cdots x_n^{m_n}$ on a
$$R_{p}(x^m)=<p,m>x^m.$$
L'entier $\delta=<p,m>=\sum_{i=1}^{n}p_im_i$ est appelé
${p}$-degré de quasi-homogénéité du monôme $x^m$. 

\begin{defi}\label{pol.qh}
Une fonction polynomiale $h$ de $\C^n$ dans \C\ est dite quasi-\-ho\-mo\-gène de poids
$p$ lorsque le $p$-degré de chacun des monômes
$x^m=x_1^{m_1}\cdots x_n^{m_n}$ de coefficient non nul de $h$ est constant. Cette
constante est le $p$-degré de quasi-homogénéité de $h$. En d'autres
termes, une fonction polynomiale est $p$-quasi-\-ho\-mo\-gène de degré $\delta$
si et seulement si $$R_{p}(h)=\delta h.$$
\end{defi}
 
Le ${p}$-ordre
$\nu_{p}$ d'une série $\sum_k a_mx^m$ est le minimum des 
${p}$-degrés des monômes de coefficient non nul de cette série. 
L'espace ${\cal P}_n$ des fonctions polynomiales sur $\C^n$ est décomposable en
somme directe des espaces ${\cal P}_n^{\delta}$ des fonctions polynomiales
${p}$-quasi-\-ho\-mo\-gènes de degré $\delta $, $\delta $ décrivant
$\N$. La notion de quasi-homogénéité s'étend aux champs de vecteurs:

\begin{defi} Un champ de vecteurs polynomial
$X=\sum_{i=1}^{n}a_i(x)\frac{\partial}{\partial x_i}$ est quasi\-ho\-mo\-gène de poids
$p$ et de degré $\delta $ si
$$[R_{p},X]=\delta X.$$
\end{defi}

Remarquons qu'avec cette convention, les champs linéaires $X$ sont quasi-\-ho\-mo\-gènes de
degré 0 et les champs constants 
$\frac{\partial}{\partial x_i}$ sont de degré négatif $-p_i$.
Dans la suite, le poids $(p_1,\cdots p_n)$ est fixé. On note $R=R_p$ et
les dénominations ''quasi-\-ho\-mo\-gène, degré, ordre'' se rapporteront toujours à ce poids.
On vérifie facilement les résultats suivants:

\begin{prop}\label{qh} 
 
i) Soit $a_k$ un polynôme quasi-\-ho\-mo\-gène de degré $k$, $X_l$ un champ
quasi-\-ho\-mo\-gène de degré $l$. Le polynôme $X_l(a_k)$ est alors quasi-\-ho\-mo\-gène, de degré
$k+l$.

ii) Soit $a_k$ un polynôme quasi-\-ho\-mo\-gène de degré $k$, $X_l$ un champ
quasi-\-ho\-mo\-gène de degré $l$. Le champ $a_kX_l$ est alors quasi-\-ho\-mo\-gène de degré
$k+l$.

iii) Soit $X_k$ un champ quasi-\-ho\-mo\-gène de degré $k$, $X_l$ un champ
quasi-\-ho\-mo\-gène de degré $l$. Le champ $[X_k,X_l]$, s'il est non nul, 
est quasi-\-ho\-mo\-gène de degré
$k+l$.

iv) Soit $h$ une fonction de deux variables quasi-\-ho\-mo\-gène de degré $\delta $. Le
champ hamiltonien $X_h=h_y\partial  _x -h_x\partial _y$
est quasi-\-ho\-mo\-gène de degré $\delta -p_1-p_2.$ 

\end{prop}

En particulier, tout champ monomial $x^m\frac{\partial}{\partial x_i}$ est
quasi-\-ho\-mo\-gène de degré
$<p,m>-p_i$. En regroupant les champs monomiaux de même degré,
on obtient une décomposition de l'espace des champs formels s'annulant à l'origine en somme
directe de sous-espaces de champs quasi-\-ho\-mo\-gènes de degré
$\delta $,
$\delta $ décrivant $\N$. Ceci permet de définir le $p$-ordre d'un tel
champ formel, ainsi que le $p$-degré d'un champ polynomial.

\bigskip

{\bf Exemple.} Le cusp $h=y^{p}-x^{q}$ est quasi-\-ho\-mo\-gène de poids 
$(p_1,q_1)$,
où $p_1$ et $q_1$ sont définis par: $p=dp_1$, $q=dq_1$, $p_1$ premier avec
$q_1$. Il est de $(p_1,q_1)$-degré $p\vee q=dp_1q_1$. Le champ 
$X_h=h_y\partial  _x -h_x\partial _y$
est $(p_1,q_1)$-quasi-\-ho\-mo\-gène de degré $dp_1q_1-p_1-q_1$.

\section{Perturbation à séparatrices fixées: champs et formes logarithmiques}

On fixe maintenant le lieu des séparatrices $S$ de $X=X_0+\cdots$
où $S$ désigne un
germe d'hypersurface analytique à l'origine de $\C^n$, d'équation réduite $h_0=0$.
On demande donc que $X$ reste tangent à $S$ c'est-à-dire que $X(h_0)$ appartienne à
l'idéal
$(h_0)$. Ceci est la définition d'un champ \textit{logarithmique}.  On désigne par
$\chi(\log S)$ le ${\cal O} _n$-module des germes de champs
 logarithmiques pour $S$.
Toutes les notions rappelées ici concernant $\chi(\log S)$ sont exposées dans
\cite{SAITO.2}. Ce module est stable par crochet de Lie. En dimension deux, il est libre de
rang deux. En dimension supérieure, ce n'est plus le cas. On a le critère suivant:

\begin{prop}\label{critere}
Les $n$ champs de vecteurs $X_i$ de $\chi(\log S)$ forment une base de ce module si
et seulement si le déterminant de leurs composantes est produit de $h_0$ par un
unité.
\end{prop}

Une $q$-forme logarithmique à pôles sur $S$ est une $q$-forme méromorphe
$\omega $ telle que $h_0\omega $ et $h_0dw$ soient holomorphes. Elles forment un 
${\cal O} _n$-module noté $\Omega ^q(\log S)$. Ces modules munis de la
différentielle extérieure usuelle forment un complexe. Les notions de produit
intérieur et dérivation de Lie s'étendent aux champs et formes
logarithmiques. En particulier, le module
$\Omega ^1(\log S)$ est dual de $\chi(\log S)$ pour l'accouplement des champs et
1-formes. Une manière duale de caractériser les hypersurfaces $S$ telles que $\chi(\log S)$ --ou 
$\Omega ^1(\log S)$)-- soit libre est la suivante:

\begin{prop}\label{flog}
Le ${\cal O} _n$-module $\Omega ^1(\log S)$ est libre si et seulement si
$$\wedge_{i=1}^{n}\Omega ^1(\log S)=\Omega ^n(\log S)$$ c'est-à-dire si et seulement
si il existe $n$ 1-formes logaritmiques $\omega _1, \cdots \omega _n$ et une
unité $u$ de ${\cal O} _n$ telles que 
$$\omega _1\wedge\cdots\wedge \omega _n=u \frac{dz_1\wedge\cdots\wedge dz_n}{h}.$$
De plus, sous cette condition, les modules $\Omega ^q(\log S)$ sont tous libres de
base  $\omega _{i_1}\wedge\cdots\wedge \omega _{i_q}.$ 
\end{prop}

Notons qu'en général, c'est-à-dire en dehors du cas où $S$ est à croisements
normaux ($\omega_i=dz_i/z_i$), les formes $\omega _i$ ne sont pas toutes fermées.

\begin{rmq}\label{derivation} Si $\chi(\log S)$ est libre de base $\{X_1,\cdots X_n\}$,
et si $\{\omega _1,\cdots \omega _n\}$ est une base de $\Omega ^1(\log S)$ duale de 
$\{X_1,\cdots X_n\}$ pour l'accouplement des champs et 1-formes, pour tout $f$ de ${\cal
O}_n$ on a
$$df=\sum_{i=1}^nX_i(f)\omega _i.$$
\end{rmq}

On se place maintenant en dimension deux. On suppose que $h_0=0$ est une
équation réduite de $S$, quasi-\-ho\-mo\-gène pour le champ $p$-radial
$R$, de degré $d_0$. On considère la ''forme-volume
logarithmique''
$$\Omega =\frac{dx\wedge dy}{h_0}.$$
Elle induit un isomorphisme $\sharp$
entre $\chi(\log S)$ et $\Omega ^{1}(\log S)$ défini par $X^{\sharp }=i_X\Omega
$ où $i_X$ désigne le produit intérieur par $X$. L'isomorphisme inverse est 
noté $\flat$. 

Soit $h_1\cdots h_p$
une décomposition de $h_0$ en éléments irréductibles et $h$ une fonction quasi-\-ho\-mo\-gène de  degré $\delta $, telle que $h=0$ soit encore 
une équation de $S$: $h=h_1^{n_1}\cdots h_p^{n_p}$. 
La forme
$$\omega _0=\frac{1}{\delta }\frac{dh}{h}$$
appartient à $\Omega ^{1}(\log S)$. En effet, 
$$h_0\omega _0=h_1\cdots h_p\frac{1}{\delta }\sum_{i=1}^{p}n_i\frac{dh_i}{h_i}$$
est holomorphe, et $h_0d\omega _0$ est trivialement holomorphe 
puisque $\omega _0$ est fermée. On s'intéresse au champ logarithmique 
$X_0$ dual de $\omega _0$.

\begin{prop}\label{X_0integrable} Soit $X_0=\omega _0^{\flat}$, et 
$X_h=h_y\partial  _x -h_x\partial _y$ le champ
hamiltonien dual de $dh$ pour la forme volume $dx\wedge dy$. 
On a

i-) $X_0=\frac{h_0}{\delta h}X_h.$ En particulier $h$ est une intégrale première du champ $X_0$.

ii-) $X_0$ est un champ à singularité isolée, quasi-\-ho\-mo\-gène pour le champ
$p$-radial
$R$, de degré $\delta _0=d_0-p$.

iii-) Les champs logarithmiques $X_0$ et $R$ forment une base de $\chi(\log S)$ et les
formes logarithmiques $-\omega _R=-R^\sharp$ et $\omega _0=X_0^{\sharp}$ forment une base
duale de $\{X_0, \ R\}$ pour l'accouplement des champs et 1-formes. En particulier, pour tout
$f$ de
${\cal O}_2$, on a
$$df=R(f)\omega _0-X_0(f)\omega _R.$$
\end{prop}

\begin{preuve} i-) De $dx\wedge dy(X_h,\cdot)=dh$, on déduit
$$\Omega (\frac{h_0}{\delta h}X_h,\cdot)=\frac{dh}{\delta h}=\omega_0,\ \mbox{d'où }\omega
_0^{\flat}=\frac{h_0}{\delta h}X_h. $$ 

ii-) se déduit de $dx\wedge dy(X_0,\cdot)=h_0\omega_0$ qui est une forme à singularité
isolée. On calcule le degré de quasi-homogénéité de $X_0$ en utilisant la proposition
(\ref{qh}).

iii-) Le critère (\ref{critere}) est vérifié d'après:
$$\Omega (X_0,R)=\frac{dx\wedge dy}{h_0}(\frac{h_0}{\delta h}X_h,R)=
\frac{dh}{\delta h}(R)=\frac{R(h)}{\delta h}=1.$$ Les autres affirmations se déduisent de
cette même égalité.
\end{preuve}

\begin{defi} Dans la décomposition unique de tout champ logarithmique en
$aX_0+bR$, nous désignerons la première composante par ''composante
intégrable'' (elle admet $h$ comme intégrale première), et
la seconde par ''composante dissipative'' (la forme duale $b\omega_R$ n'est
jamais fermée). Nous dirons de plus que la composante intégrable est
''hamiltonienne'' lorsqu'elle est duale d'une forme logarithmique de type
$dg/g$ (c'est par exemple le cas lorsque $a=a(h)$).
\end{defi}

\bigskip

{\bf Exemple 1.} (Poincaré-Dulac):

On prend: $h_0=xy$, et $h=x^py^q$ où $p$ et $q$ sont des entiers strictement
positifs premiers entre eux. Les fonctions $h_0$ et
$h$ sont
évidemment homogènes mais aussi quasi-\-ho\-mo\-gènes pour le champ $(q,p)$-radial
$R=qx\partial _x+py\partial _y$. Pour ce champ, on a $d_0=\deg h_0=p+q$,
$\delta\deg h =2pq$,
et $\delta _0=\deg X_0=0$.
Le champ hamiltonien $X_h$ n'est pas à singularité isolée
dès que $p$ ou $q$ est strictement supérieur à 1. On a ici:
$$\omega _0=\frac{1}{2pq}(p\frac{dx}{x}+q\frac{dy}{y})\ \ \mbox{et 
}X_0=\frac{1}{2pq} (qx\partial_x-py\partial_y).$$
La base duale de $X_0$ et $R$ est formée de 
$$-\omega _R=p\frac{dx}{x}-q\frac{dy}{y}\ \ \mbox{et de }\omega _0.$$

\bigskip

{\bf Exemple 2.} (Les singularités de type noeud-col):
On pose ici $X_0=y\partial_y$. Le champ
$X_0$ est quasi-\-ho\-mo\-gène pour le champ $(1,0)$-radial $R=x\partial _x$
de degré $\delta _0=0:$
$[x\partial _x, y\partial _y]=0$. 
Deux choix sont ici possible pour la séparatrice $S$: $x=0$ ou $xy=0$.
Dans le premier cas, les champs $X_0$ et $R$ ne forment pas une base de
$\chi(\log S)$. Le contexte où nous nous sommes placés ici impose donc de ne
considérer que les singularités de type noeud-col qui conservent les deux axes
comme séparatrices.

\bigskip

Lorsque l'intégrale première $h$ est elle-même à singularité isolée, on a
$h_0=h$, et $X_0=\delta^{-1}X_h$. En particulier:

\bigskip

{\bf Exemple 3.} (Le cusp; voir \cite{LOR} et \cite{ZOL}):

On a ici: $h_0=h=y^{p}-x^{q}$ où $p$ et $q$ sont des entiers strictement
positifs. Dans ce cas, on a vu que $h$ est quasi-homogène pour le champ radial
$R=p_1x\partial _x+q_1y\partial _y$, ($p=dp_1$, $q=dq_1$, $p_1$ premier avec
$q_1$), avec $d_0=\delta =p\vee q=dp_1q_1$ et
$\delta _0=dp_1q_1-p_1-q_1$.
La base duale de $X_0$ et $R$ est formée de $-\omega _R$ et $\omega_0$ avec
$$\omega _R=\frac{1}{h}(-q_1ydx+p_1xdy)\ \ \mbox{et  }
\omega_0=\frac{1}{ \delta}\frac{dh}{ h}.$$

\bigskip

On notera $\widehat\chi(\log S)$, $\widehat\Omega ^1(\log S)$ les complétés
formels de $\chi(\log S)$ et 
$\Omega ^1(\log S)$. Les résultats
ci-dessus s'étendent dans ce contexte formel.

\section {Pré\-nor\-ma\-li\-sation des champs et feuilletages}

On se donne un sous-module $M$ du module sur $\widehat{\mathcal O}_n$ des champs de vecteurs
formels, gradué par un degré de quasihomogénéité $k$: $M=\oplus _k M_k$. On note
$M_{\leq k}=\oplus_{l\leq k} M_k$ le module des champs polynomiaux de degré au plus $k$.
Soit $X=X_0+\cdots$ un champ dont la partie initiale $X_0$ est
quasi-homogène de degré $\delta_0
$. On considère le complexe de Koszul:
$$0\longrightarrow M \stackrel{d^0}{\longrightarrow }
M \stackrel{d^1}{\longrightarrow }
\wedge^2 M \longrightarrow\cdots$$
L'opérateur $d^0$ est $\mbox{ad}_{X_0}=[X_0,\cdot]$. L'opérateur $d^1$ est
ici nul (on ne considère qu'un seul champ de vecteurs) et rend l'appellation
''complexe de Koszul'' un peu superflue. Nous la conservons néanmoins car
devient nécessaire dès qu'on cherche à normaliser une famille de plusieurs champs de
vecteurs (voir
 \cite{STOLO}). Dans les notations suivantes, l'indice $K$ rappelle l'utilisation de ce complexe
de Koszul:
\begin{eqnarray*}
B_K^1&=& \mbox{Im}(d^0)\\
Z_K^1&=& \mbox{Ker}(d^1)\ \ \ (= M \mbox{ ici})\\ 
H_K^1&=&  Z_K^1/ B_K^1.
\end{eqnarray*}
Le théorème suivant généralise un fait exposé par J. Martinet au
séminaire Bourbaki \cite{MAR} dans le cas d'un champ initial linéaire semi-simple. 

\begin{lem}[pré\-nor\-ma\-li\-sation des champs de vecteurs]\label{Martinet.1} On se donne un
sup\-plé\-men\-taire
$V$ de $B_K^1$ dans $M$. Il existe un 
champ $Y$ appartenant à $V$ et un difféo\-mor\-phisme formel $\Phi $ tel que
$$\Phi ^*X=X_0+Y.$$
\end{lem}

\begin{preuve}
Le problème étant formel, on cherche une conjugaison $\Phi $ 
sous forme $\exp Z$ où $Z$ est un champ de vecteurs formel. On a
$$(\exp Z)^*X=\sum_{k\geq0}\frac{\mbox{ad}_Z^{(k)}}{k!}(X)=X+[Z,X]+\cdots.$$
Remarquons d'abord que d'après la proposition (\ref{qh}), l'opérateur
$d^0=\mbox{ad}_{X_0}$ est compatible avec la décomposition par le degré, avec un
décalage de $\delta_0=\mbox{deg}(X_0)$. Ceci nous permet de décomposer par le
degré les espaces $B_K^1$, $H_K^1$, et $V$. 
Dans ce qui suit, les indices correspondent au degré du champ correspondant.
L'existence de $\Phi $ s'établit par la récurrence suivante. 
On suppose qu'on a une
telle décomposition jusqu'en degré $k$, c'est-à-dire qu'il existe un champ
$Y_k$ appartenant à $V_{\leq k}$ et un champ
$Z_{k-\delta_0}$ dans $M_{\leq k-\delta_0}$ tels que
$$(\exp Z_{k-\delta_0})^*X=X_0+Y_k+\widetilde{Y}_{k+1}+\cdots $$
où $\widetilde{Y}_{k+1}$ appartient à  $M_{k+1}$.
On décompose ce premier terme du reste dans $M_{k+1}=Im(d^0)_{k+1} \oplus V_{k+1}$:
$$\widetilde{Y}_{k+1}=d^0(Z_{k+1-\delta_0 })+Y_{k+1}.$$
En composant à gauche la conjugaison précédente par
$\exp(Z_{k+1-\delta_0})$, on fait alors disparaitre le terme
$d^0(Z_{k+1-\delta_0 })$. En effet, en notant $\equiv$ l'égalité modulo des
termes d'ordre supérieur ou égal à $k+2$, on a
\begin{eqnarray*}
\exp(Z_{k+1-\delta_0})^*\exp(Z_{k-\delta_0 })^*X&\equiv &\exp(Z_{k-\delta_0})^*X+
[Z_{k+1-\delta_0 },\exp(Z_{k-\delta_0})^*X]\\
&\equiv &X_0+Y_k+Y_{k+1}+d^0(Z_{k+1-\delta_0 })+[Z_{k+1-\delta_0},X_0]\\
&\equiv &X_0+Y_k+Y_{k+1}
\end{eqnarray*}
ce qui démontre le lemme.
\end{preuve}

\begin{rmq}\label{unicite.1}
i) Il n'y a pas d'énoncé analytique analogue à ce lemme, car nous avons utilisé de
manière essentielle le fait qu'une conjugaison formelle est 
l'exponentielle d'un champ formel, fait que nous perdons dans le contexte analytique. 

ii) Pour un supplémentaire $V$ fixé, ce
processus ne donne pas un unique champ équivalent à $X$ dans $V$. En effet,
à chaque étape, on peut rajouter au champ $Z_{k+1-\delta_0 }$ un champ commutant avec la partie
initiale $X_0$. Ce choix perturbe les termes d'ordre supérieur et modifie donc la forme
normale formelle obtenue dans $V$. Pour cette raison, nous parlons seulement
ici de ''prénormalisation''.
\end{rmq}

Considérons maintenant la relation d'équivalence formelle 
classifiant les champs \textit{à unité près}:
$$X\approx X' \Leftrightarrow \exists \Phi \in \widehat{\mbox{Diff}}(\C^2,0),\ \exists
u \in \widehat{\cal O}_2,\ u(0)\neq 0,\ \Phi ^*X=uX'.$$

\begin{lem}[pré\-nor\-ma\-li\-sation des feuilletages]\label{Martinet.2} On se
donne un sup\-plé\-men\-taire
$W$ de $A_K^1=B_K^1+\widehat{\mathcal O}_nX_0$ dans $M$. 
Il existe un 
champ $Y$ appartenant à $W$, un difféo\-mor\-phisme formel $\Phi $ et une
unité formelle $u$ tels que
$$\Phi ^*X=u(X_0+Y).$$
\end{lem}

\textit{Preuve.}
On reprend le raisonnement par récurrence de (\ref{Martinet.1}) avec l'hypothèse suivante:
on suppose qu'il existe un champ
$Y_k$ appartenant à $W_{\leq k}$, un champ
$Z_{k-\delta_0}$ dans $M_{\leq k-\delta_0}$ et un polynôme $u_{k-\delta_0 }=1+\cdots$ de degré
inférieur ou égal à $k-\delta_0 $ tels que
$$(\exp Z_{k-\delta_0})^*X=u_{k-\delta_0 }(X_0+Y_k)+\widetilde{Y}_{k+1}+\cdots $$
où $\widetilde{Y}_{k+1}$ appartient à  $M_{k+1}$.
On décompose ce premier terme du reste en
$$\widetilde{Y}_{k+1}=d^0(Z_{k+1-\delta_0 })+a_{k+1-\delta_0 }X_0+Y_{k+1}$$
où le champ $Y_{k+1}$ appartient à $W_{k+1}$ et $a_{k+1-\delta_0 }$ est quasi-homogène de
degré ${k+1-\delta_0 }$. La conjugaison de $X$ par
$\exp(Z_{k+1-\delta_0})\circ \exp (Z_{k-\delta_0 })$ donne, modulo des
termes d'ordre supérieur ou égal à $k+2$:
\begin{eqnarray*}
\exp(Z_{k+1-\delta_0})^*\exp(Z_{k-\delta_0 })^*X&\equiv & 
u_{k-\delta_0 }(X_0+Y_k)+d^0(Z_{k+1-\delta_0 })+Y_{k+1}+a_{k+1-\delta_0}X_0+\\
&&+[Z_{k+1-\delta_0},(1+\cdots)X_0]\\
&\equiv & u_{k-\delta_0 }(X_0+Y_k)+Y_{k+1}+a_{k+1-\delta_0}X_0\\
&\equiv & (u_{k-\delta_0 }+a_{k+1-\delta_0})(X_0+\frac{u_{k-\delta_0 }}{u_{k-\delta_0 }+a_{k+1-\delta_0}}Y_k
+\\
&&+\frac{1}{u_{k-\delta_0 }+a_{k+1-\delta_0}}Y_{k+1})\\
&\equiv & (u_{k-\delta_0 }+a_{k+1-\delta_0})(X_0+Y_k+Y_{k+1})
\end{eqnarray*}
la dernière égalité provenant de
$$\frac{u_{k-\delta_0 }}{u_{k-\delta_0 }+a_{k+1-\delta_0}}Y_k=(1-a_{k+1-\delta_0}+\cdots)Y_k\equiv Y_k$$
puisque $2k+1-\delta_0 >k+1$ dès que $k>\delta_0 $, et de
$$\frac{1}{u_{k-\delta_0 }+a_{k+1-\delta_0}}Y_{k+1}=(1+\cdots)Y_{k+1}\equiv
Y_{k+1}.
\ \carre $$

\begin{rmq}\label{unicite.2}
Là encore, il n'y a pas unicité de la conjuguante et de la forme
normale formelle du feuilletage ainsi obtenues. 
Nous avons même plus de liberté dans ces choix que dans la situation précédente puisque,
à chaque étape, on peut rajouter au champ
$Z_{k+1-\delta_0 }$ un champ qui ne commute avec la partie initiale $X_0$ que modulo un
champ multiple de $X_0$.
\end{rmq}

\bigskip

Nous nous plaçons maintenant en dimension 2, et considérons une perturbation
$X$ du champ quasi-homogène $X_0$ de degré $\delta_0$ préservant un ensemble
analytique invariant $S$. Le champ $X$ appartient donc au module
$M=\widehat{\chi}(\log S)$ des champs logarithmiques pour $S$. Nous supposons que les champs $X_0$ et $R$ forment une base de
$\widehat{\chi}(\log S)$: c'est par exemple le cas lorsque le champ logarithmique $X_0$ est dual d'une 1-forme
$dh/h$, d'après la proposition (\ref{X_0integrable}).
Le champ $X$ s'écrit
$$X=aX_0+bR, $$
avec $a(0)=1$ et $\nu(bR)=\nu(b)+0>\delta_0$ où $\nu$ désigne l'ordre de
quasi-homogénéité relatif à $R$.
Dans ce contexte, les lemmes de pré\-nor\-ma\-lisation (\ref{Martinet.1}) et
(\ref{Martinet.2}) peuvent être précisés de la manière suivante.
Remarquons d'abord que, puisque $W$ est un supplémentaire d'un espace contenant
$\widehat{\cal O}_2\cdot X_0$, $W$ n'a pas de  composante intégrable, et est donc inclus dans 
$\widehat{\cal O}_2\cdot R$: \textit{la composante dissipative (en $R$) d'un
  espace de formes prénormales du champ $X$ est une espace de formes
  prénormales du feuilletage défini par $X$.} 

Nous pouvons de plus nous limiter à ne considérer
que des conjugaisons
du type suivant:

\begin{defi} Nous dirons qu'un difféomorphisme formel est fibré par rapport au champ $R$
s'il est l'exponentielle d'un champ multiple de $R$.
\end{defi}

Un tel difféomorphisme préserve donc le feuilletage défini par le champ $R$,
feuilletage qui est transverse au feuilletage initial défini par le champ
$X_0$ en dehors de $S$.
%Dans la suite, un difféomorphisme fibré (sans autre précision) désignera un
%difféomorphisme fibré par rapport au champ radial $R$.

\begin{lem}
[pré\-nor\-ma\-li\-sation 
à séparatrices 
fixées]
\label{Martinet.3} 

On se donne un sup\-plé\-men\-taire $V$ de $B^1_K$ dans $M=\widehat{\chi}(\log
S)$. Soit $W\subset \widehat{\mathcal O}_2R$ la composante dissipative de $V$,
supplémentaire de $A_K^1=B_K^1+\widehat{\mathcal O}_2X_0$ dans
$M$.  

1- Il existe un champ $Y$ appartenant à $V$, un difféo\-mor\-phisme formel
$\Phi$, tels que $\Phi ^*X=X_0+Y.$

2- Il existe un 
champ $Y$ appartenant à $W$, un difféo\-mor\-phisme formel $\Phi$ fibré par
rapport à $R$ et une unité formelle $u$ tels que
$\Phi ^*X=u(X_0+Y).$
\end{lem}

\begin{preuve} Le module $M=\widehat{\chi}(\log S)$ étant stable par crochet
 de Lie, cet énoncé est un corollaire des deux lemmes précédents
 \ref{Martinet.1} et  \ref{Martinet.2}, sous réserve de vérifier
de plus que la conjugaison des feuilletages peut être choisie de manière
fibrée. On modifie l'argument de récurrence de  \ref{Martinet.2} en supposant que $\exp Z_{k-\delta_0 }$ est fibrée, de la
forme $\exp b_{k-\delta_0 }R$, où $b_{k-\delta_0 }$ est un polynôme de degré au plus
$k-\delta_0 $. Dans la décomposition du premier terme du reste
$$\widetilde{Y}_{k+1}=d^0(Z_{k+1-\delta_0 })+a_{k+1-\delta_0 }X_0+Y_{k+1}$$
l'écriture de la somme des deux premiers termes n'est pas unique. On peut ajouter à $Z_{k+1-\delta_0 }$
une composante intégrable quasi-homogène de degré 
$k+1-\delta_0 $ arbitraire: cela n'affectera que la valeur du coefficient $a_{k+1-\delta_0 }$. Nous
utilisons cette liberté pour annuler la composante intégrable de $Z_{k+1-\delta_0 }$ et
choisir ce champ sous la forme $b_{k+1-\delta_0 }R$.
\end{preuve}

\begin{rmq}\label{unicite.3} Deux champs fibrés
par rapport à $R$ ont même image par $d^0$ si et seulement si ils diffèrent d'un champ
fibré $cR$ tel que $[X_0,cR]=0$, c'est-à-dire d'un champ fibré $cR$ vérifiant
$X_0(c)=0$. On peut donc faire opérer sur l'ensemble des formes
prénormales d'un feuilletage donné le
groupe des conjugaisons fibrées de la forme
$\exp cR$
où $c$ est une intégrale première de la partie initiale $X_0$.
\end{rmq}

\bigskip

D'après le lemme de prénormalisation \ref{Martinet.3}, les  quotients qui
nous intéressent sont maintenant
$$H_K^1=\widehat{\chi}(\log S)/B_K^1$$
pour la classification des champs à séparatrices fixées,  et
$$G_K^1=\widehat{\chi}(\log S)/A_K^1\ \mbox{ avec  } A_K^1=\widehat{\cal O}_2\cdot X_0+B_K^1$$
pour la classification de
feuilletages à séparatrices fixées. Remarquons que ce sont des modules sur
l'anneau ${\mathcal I}$ des
intégrales premières de $X_0$.
On désigne par $\mbox{Ker}(X_0)$ et $\mbox{Coker}(X_0) =\widehat{\cal
O}_2/\mbox{Im}(X_0)$ les noyau et conoyau du champ $X_0$ vu comme
dérivation agissant sur $\widehat{\cal O}_2$. 
Ce sont aussi des ${\mathcal I}$-modules.
%On note $\mbox{Coker}(X_0)_\delta$ le sous-module des éléments $b$
%de $\mbox{Coker}(X_0)$ d'ordre strictement supérieur à l'entier $\delta$.

\begin{thm}\label{PFNF}
Soit $X_0$ un champ quasi-homogène de degré $\delta_0$ par rapport à un champ
quasi-radial $R$, laissant invariant une courbe $S$. On suppose que les champs
$X_0$ et $R$ forment une base des champs logarithmiques pour $S$.
\begin{enumerate}
\item 
%Le module des formes prénormales de feuilletages définis par
%  $X_0+\cdots$ est isomorphe au module 
$G_K^1=\mbox{Coker}(X_0)\cdot R$.
\item Si le degré $\delta_ 0$ du champ initial est nul (exemples 1 et 2), 
%le module des formes prénormales du champ  $X_0+\cdots$ est
$$H_K^1=\mbox{Coker}(X_0)\cdot X_0+\mbox{Coker}(X_0)\cdot R.$$
\end{enumerate}
\end{thm}

\begin{preuve}
De la relation de quasi-homogénéité $[R,X_0]=\delta_0X_0$, on déduit
$$[X_0,aX_0+bR]=(X_0(a)-\delta_0b)X_0+X_0(b)R.$$
Un champ $X=\alpha X_0+\beta R$ appartient donc à $B^1_K$ si et seulement si
il existe deux coefficients $a$ et $b$ solutions du système
$$\left\{\begin{array}{lll}
X_0(a)-\delta_0b &=&\alpha  \\
X_0(b)& =& \beta 
\end{array}\right.$$
En particulier, un champ $X=\alpha X_0+\beta R$ appartient à 
$A^1_K=\widehat{\cal O}_2\cdot X_0+B_K^1$ si et seulement si
il existe un  coefficient $b$ solution de
$$X_0(b) = 
\beta .$$
L'application
\begin{eqnarray*}
\widehat{\chi}(\log S)\ \ \ \ & \rightarrow& \widehat{\cal O}_2/\mbox{Im}(X_0)\cdot R\\
X=\alpha X_0+\beta R &\mapsto& [\beta ]\cdot R
\end{eqnarray*}
a pour noyau $A_K^1$ et définit donc un isomorphisme entre $G_K^1$ et
$\mbox{Coker}(X_0)\cdot R$.

Si le degré $\delta_0$ est nul, le système précédent est diagonal, et on
constate immé\-dia\-tement que $X=\alpha X_0+\beta R$ appartient à $B^1_K$ si et
seulement si il existe deux coefficients $a$ et $b$ solutions de
$X_0(a)=\alpha$ et de $X_0(b)=\beta$. On identifie donc comme ci-dessus
$H_K^1$ et $\mbox{Coker}(X_0)\cdot X_0+\mbox{Coker}(X_0)\cdot R$.
\end{preuve}

\bigskip

Remarquons que dès que le degré $\delta_0$ du champ initial est strictement
positif, le problème de la détermination des formes prénormales de champs
devient plus difficile: il ne se décompose plus en somme directe des deux
problèmes ''détermination d'une classe de feuilletages'' puis ''classification
des champs au sein de cette classe''.

Le théorème \ref{PFNF} ramène donc la détermination des espaces de formes
pré\-nor\-males formelles de feuilletages (et de champs lorsque $\delta_0=0$) au
calcul de $\mbox{Coker}(X_0)$.
Sur les exemples ''classiques'', le conoyau de $X_0$ se calcule aisément.
Reprenons les exemples du paragraphe 2:

\bigskip

{\bf Exemple 1.} (Poincaré-Dulac):

Soit $X_0=\frac{1}{2pq} (qx\partial_x-py\partial_y)$, quasi-homogène pour 
$R=qx\partial _x+py\partial _y$. Les seules
obstructions à la résolution de l'équation $X_0(f)=g$ proviennent des termes
$g_{kp,kq}$ de la série de Taylor de $g$. Le conoyau de $X_0$ est donc ici
$\C[[h]]$, avec $h=x^py^q$. Nous en déduisons donc qu'il existe des séries
formelles d'une variable $a$ et $b$ avec $\nu(b)>0$ telles que

 $$X\sim a(x^py^q)X_0+b(x^py^q)R,  \ \ X\approx X_0+b(x^py^q)R.$$

\bigskip

{\bf Exemple 2.} (Les singularités de type noeud-col):

On a ici: $X_0=y\partial_y$ et $R=x\partial _x$.
On vérifie immédiatement que le conoyau est ici $\C[[x]]$. Donc pour tout
champ $X=X_0+\cdots$  il existe des séries
formelles d'une variable $a$ et $b$ avec $\nu(b)>0$ telles que
 $$X\sim a(x)X_0+b(x)R,  \ \ X\approx X_0+b(x)R.$$

{\bf Exemple 3.} (Le cusp $h=y^p-x^q$):

Soit $X_0= \frac{1}{ \delta}(py^{p-1}\partial_x+qx^{q-1}\partial_y)$.
Le conoyau de cette dérivation a été calculé par F. Loray dans \cite{LOR}: il
est engendré sur $\C[[h]]$ par les monômes $a_{k,l}=x^ky^l$, $k=0, \cdots
q-2$, $l=0\cdots p-2$. L'espace des formes
prénormales de feuilletages obtenues ici est donc
%\sum_{\begin{array}{c}{\scriptstyle k=0,\cdots q-2,}\\ {\scriptstyle l=0\cdots p-2}\end{array}}
$$X\approx X_0 +\sum_{k=0,\cdots q-2,\ l=0\cdots p-2}
 \C[[h]]x^ky^l\ R.$$

\bigskip

Nous nous proposons dans le paragraphe suivant, de généraliser ce dernier exemple à toute dérivation
$X_0=X_h$, où $h$ est une fonction quasi-homogène à singularité isolée.

\section{Conoyau d'une dérivation hamiltonienne $X_0$}

Nous nous plaçons sous l'hypothèse (2-) mentionnée dans l'introduction: le champ initial 
$X_0$ est hamiltonien, dual pour la forme $dx\wedge dy$ d'une forme
$\delta^{-1}dh$, où l'intégrale première $h$ est quasi-homogène de degré
$\delta$. Nous savons d'après (\ref{X_0integrable}) que les champs $X_0$ et $R$ forment une base
de $\widehat\chi(\log S)$.
Nous supposons de plus que $h$ est à singularité isolée.
On désigne par $J_0$
l'idéal des composantes de $X_0$ dans la base usuelle (qui est aussi l'idéal jacobien de $h$), et par $a_1=1,\cdots
a_{\mu}$, une base engendrant l'espace vectoriel ${\cal O}^2/J_0$.

Dans ce paragraphe, nous étudions le champ $X_0$ vu comme opérateur de
dérivation sur l'anneau de fonctions $\C\{h\}$ (ou $\C[[h]]$). En d'autres
termes, nous écrivons ici $X_0$ pour désigner la dérivée de Lie ${\mathcal L}_{X_0}$.

\begin{thm}\label{conoyau}\footnote{R. Schaefke a communiqué à l'auteur une
autre preuve de ce résultat, reposant néanmoins sur un lemme analogue au lemme
de division \ref{lemme de division}.}
Le noyau $\mbox{Ker}(X_0)$ de la
dérivation $X_0$ est l'anneau ${\mathcal I}=\C\{h\}$ des intégrales premières de $X_0$. 
Son conoyau ${\cal O}_2/\mbox{Im}(X_0)$ est le module libre de rang $\mu$ engendré sur 
$\C\{h\}$ par
$(a_1, a_2,\cdots a_{\mu})$.

On a un résultat analogue pour la dérivation $X_0$ agissant sur l'anneau des séries
formelles
$\widehat{\cal O}_2$, en substituant dans l'énoncé ci-dessus l'anneau $\C\{h\}$ par
$\C[[h]]$.
\end{thm}

\begin{preuve}
Clairement $\C\{h\}$ est inclus dans le noyau
de $X_0$. L'égalité résulte d'un théorème de Mattei-Moussu (\cite{MM}):
Puisque $h$ est à singularité isolée, elle 
n'est pas une puissance d'une autre intégrale première et toute intégrale première analytique (resp. formelle) du
champ est obtenue par composition à gauche de $h$ par une série convergente (resp.
formelle).

Pour déterminer le conoyau de $X_0$, nous nous inspirons d'un argument de B. Malgrange utilisé
dans le cadre des modules de cohomologie relative à un germe d'application holomorphe:
\cite{MAL}. On considère ici les
$\C\{h\}$-modules
$$E=
%\frac{J_0}{ \mbox{Im}(X_0)}=
\frac{J_0}{\mbox{Im}(X_0)} \    \ \mbox{et  } 
F=\frac{{\cal O}_2}{\mbox{Im}(X_0)}.$$

Le quotient $F/E$ est un
espace vectoriel de dimension finie
$\mu$. Ces modules sont munis d'une E-F connexion $\nabla :E    \rightarrow F$ (au sens défini dans
\cite{MAL}) que l'on construit à l'aide d'un lemme de division:

\begin{lem}[lemme de division]\label{lemme de division}
Pour tout élément $f$ de $J_0$, il existe $a$ et $b$ dans ${\cal O}_2$ tels que
$$f=ah+X_0(b).$$
De plus, dans cette écriture, le coefficient $a$ représente un unique élément de $F$.
\end{lem}

\textit{Preuve du lemme de division.} 
On remarque d'abord que pour toute
constante $c$ non entière négative, les opérateurs $R+c$ sont bijectifs. En effet, 
en utilisant les
décompositions en composantes quasi-\-ho\-mo\-gènes $\alpha =\sum \alpha _i$, $\beta =\sum \beta _i$,
l'égalité $(R+c)(\alpha )=\beta $ équivaut à $(i+c)\alpha _i=\beta _i$ pour
tout $i$, système que
l'on résoud formellement ou analytiquement sans autre obstruction que celle annoncée.
De plus, pour toute constante
$c$ on a la relation de commutation
\begin{eqnarray}
(R+c)X_0=[R,X_0]+X_0R +cX_0=\delta _0X_0+X_0R +cX_0=X_0(R+c+\delta _0).
\end{eqnarray}
ou encore, lorsque $R+c$ est inversible,
\begin{eqnarray}\label{com}
(R+c)^{-1}X_0=X_0(R+c+\delta _0)^{-1}.
\end{eqnarray}
Pour tout élément $f$ de $J_0$, il existe un champ $X$ tel que $f=X(h)$. Choisissons un coefficient $\alpha $ de sorte que le champ $X-\alpha R$ soit de
divergence nulle. Ceci est possible car, si $p$ désigne la somme des poids de $R$,
\begin{eqnarray*}
\mbox{div}(X-\alpha R)=0 &\Leftrightarrow& \mbox{div}(X)-R(\alpha )-\alpha \mbox{div}(R)=0\\
&\Leftrightarrow& (R+p)(\alpha )=\mbox{div}(X).
\end{eqnarray*}
et on utilise la surjectivité de l'opérateur $R+p$. Pour tout
$f$ de $J_0$ on a
$$f=X(h)= \alpha R(h) + (X-\alpha R)(h)= ah +Y(h).$$
Le champ $Y$ étant de divergence nulle, il existe
une fonction $\tilde b$ telle que 
$$Y=\tilde b_y\partial _x- \tilde b_x\partial _y.$$
On a donc $Y(h)=\tilde b_yh_x-\tilde b_xh_y=X_h(\tilde b)=X_0(b),$ pour
$b=\delta\tilde b$, d'où l'existence de la décomposition annoncée.
Pour établir son unicité, il nous faut démontrer que si $ah=-X_0(b),$
alors $a$ est dans l'image de $X_0$. L'égalité précédente s'écrit encore
$$(a    \frac{R}{\delta }+Y)(h)=0$$
et $X_0$ étant à singularité isolée, nous en déduisons l'existence de $c$ dans
${\cal O}_2$ tel que
$$a     \frac{R}{\delta }+Y=cX_0.$$
Puisque $X_0$ et $Y$ sont de divergence nulle, on obtient en appliquant l'opérateur de divergence,  
$$\frac{R+p}{\delta }(a)=X_0(c)$$
d'où, en utilisant la relation de commutation (\ref{com}),
$$a=    (R+p)^{-1}X_0(\delta  c)=X_0(R+p+\delta_0)^{-1}(\delta  c).    \ \ \carre $$

\bigskip

\noindent \textit{Fin de la preuve du théorème (\ref{conoyau}).} Ce lemme de division nous
permet de définir l'opérateur 
\begin{eqnarray*}
\nabla:\ \frac{J_0}{ \mbox{Im}(X_0)}&\longrightarrow &\frac{{\cal O}_2}{
\mbox{Im}(X_0)}\\
\end{eqnarray*}
par la formule $\nabla = \frac{R}{\delta  h}$. En effet, pour tout élément de la source
représenté par $f=ah$,
$$\nabla (f)= \frac{R}{\delta  h}(f)=\frac{R}{\delta h}(ah)=\frac{R+\delta  }{\delta }(a),$$ 
est bien défini
dans le quotient $F=\frac{{\cal O}_2}{\mbox{Im}(X_0)}.$ Cet opérateur $\C$-linéaire est une
connexion de
$\C\{h\}$-modules:
$$ \nabla(l(h)f)=l'(h)\frac{R}{\delta  h}(h)f+l(h)\frac{R}{\delta  h}(f)=l'(h)f+l(h)\nabla(f).$$
De plus, l'opérateur $\nabla$ est un isomorphisme: ceci résulte de la bijectivité de
l'opérateur $R+\delta $ et de l'unicité du coefficient $a$ obtenu au lemme de division.
L'indice
$\chi (\nabla,{\cal O}_2)$ de cet opérateur est donc nul. Un théorème d'indice analytique
pour les $E$-$F$-connexions (voir 
\cite{MAL}, théorème (3.2) page 408) nous assure alors que le rang $r$ du module $E$ est
donné par
$$r=\chi (\nabla,{\cal O}_2)+\mbox{dim}_{\C}F/E=\mu.$$

Pour trouver un système générateur du module $F$, il suffit de remarquer que le
module $E$ est isomorphe à ${\cal M}\cdot F$ où ${\cal M}$ désigne l'idéal
maximal de $\C\{h\}$. En effet, le lemme de division nous permet de définir une
application injective
\begin{eqnarray*}
E &\longrightarrow & {\cal M}\cdot F\\
f&\mapsto & ah.
\end{eqnarray*}
La surjectivité de cette application provient de la quasi-homogénéité de $h$.
L'espace vectoriel $F/{\cal M}\cdot F$ est donc isomorphe à $F/E={\cal
O}_2/J_0$ et on obtient un système générateur de $F$ à partir de
représentants d'une base de cet espace vectoriel par le lemme de Nakayama.

Enfin, il nous reste à vérifier que $F$ est sans torsion. Soit $f$ un
représentant d'un élément de
$F$ tel que $hf$ soit nul dans $F$. Puisque $hf$ appartient à ${\mathcal M}\cdot F=E$, on peut
lui appliquer l'opérateur $\nabla=R/\delta  h$, et on a $\nabla(hf)=0$ dans $F$, ce qui
signifie qu'il existe un élément $g$ de ${\mathcal O}_2$ tel que 
$$\nabla(hf)=X_0(g).$$
Les égalités intermédiaires du calcul qui suit ont un sens après avoir tensorisé par le
corps des fractions de
$\C\{h\}$:
$$\nabla (hf)=f+h\nabla f=f+h\frac{R}{\delta h}f=\frac{R+\delta }{\delta }(f).$$
On obtient donc en utilisant la relation (\ref{com}),
$$f=(R+\delta )^{-1}X_0(\delta g)=X_0(R+\delta +\delta _0)^{-1}(\delta g)$$
d'où $f=0$ dans ${\mathcal O}_2/J_0$.

Pour obtenir une version formelle de ce résultat, nous devons nous assurer que l'indice formel
$\chi (\nabla,\widehat{\cal O}_2)$ reste égal à l'indice analytique. D'après \cite{MAL.2},
ceci équivaut à vérifier que l'origine est une singularité régulière pour
l'opérateur $\nabla$. Pour cela, nous choisissons une base de $E$ engendrée par des $a_i$ qui
sont des \textit{monômes}: $a_i=x^{k_i}y^{l_i}.$ Ceci est toujours possible en utilisant la
notion de base standard considérée par J. Briançon et A. Galligo dans
\cite{BG}. On peut supposer de plus que $a_1=1.$ Calculons
le système différentiel associé à $\nabla$ dans cette base, sur le corps des fractions de
$\C\{h\}$. On a
$$\nabla(a_i)=\frac{R}{\delta  h}(a_i)=r _i\frac{a_i}{h}$$
où $r _i$ est le rationnel positif $\frac{p_1k_i+p_2l_i}{\delta }$. On a donc
%\begin{eqnarray*}\label{calcul.de.nabla}
$$\nabla(\sum_{i=1}^{\mu}d_i(h)a_i)=\sum_{i=1}^{\mu}(\frac{d}{dh}d_i(h)
a_i+d_i(h)r _i\frac{a_i}{h})$$ d'où le système diagonal 
$$h\frac{d}{dh} d_i(h)+r_id_i(h)=0,\ \ i=1\cdots \mu $$
qui est à point singulier régulier.
\end{preuve}

\bigskip

\begin{rmq}\label{solutions} Les solutions de ce système sont de la forme 
$$ d_i(h)=c_ih^{-r_i},\ \ c_i\in \C.$$
\end{rmq}

Le théorème (\ref{PFNF}) et la détermination du conoyau de $X_0$
(\ref{conoyau}) démontrent le théorème 2 énoncé dans l'introduction.

\section{Formes prénormales et cohomologie relative à la forme initiale duale}

Soit $\omega _0=\delta^{-1}dh/h$ la forme logaritmique fermée duale de
$X_0$. Le complexe logarithmique formel relatif à $\omega _0$ est défini par:
$$\widehat\Omega ^{\cdot}_{DR}(\log S):= 
\frac{\widehat\Omega ^{\cdot}(\log S)}{\omega _0\wedge \widehat\Omega
^{\cdot -1}(\log S)}$$ muni de la différentielle $d$ usuelle
(celle-ci passe au quotient puisque
$\omega _0$ est fermée). L'indice ''$\scriptstyle DR$'' pour De Rham Relatif fait référence
à cette différentielle. La cohomologie logarithmique relative à $\omega
_0$ est la cohomologie de ce complexe. Remarquons que, contrairement au complexe relatif usuel (non
logarithmique) toute 1-forme est ici relativement fermée. En effet, puisque $\omega _0\wedge \omega _{R}$ est une base de $\Omega
^{2}(\log S)$ (voir \ref{flog}), toute 2-forme logarithmique est multiple de $\omega _0$, et on a
$\widehat\Omega ^{2}_{DR}(\log S)=0.$ On a donc:
$$\widehat H^1 _{DR}(\log S)= 
\frac{\widehat\Omega ^{1}(\log S)}{\widehat B^{1}_{DR}(\log S)}$$ avec
$$\widehat B^{1}_{DR}(\log S)=\{\alpha \in \widehat\Omega ^{1}(\log
S), \ \exists f,  \ g\in \widehat{\cal O}_2,\  \alpha =df+g\cdot \omega _0 \}.$$On remarquera que $\widehat B^{1}_{DR}(\log S)$ et $\widehat
H^{1}_{DR}(\log S)$ sont encore des
$\C[[h]]$-modules.

\begin{prop}\label{calcul.Hrel} Le module
$\widehat H^{1}_{DR}(\log S)$ est le $\C [[h]]$-module dual du module des
formes prénormales de feuilletages $G_K^1$: 
$\widehat H^{1}_{DR}(\log S)=\mbox{Coker}(X_0)\cdot \omega_R$.
En particulier, il est libre de rang
$\mu$, engendré par les formes $a_i\omega _{R}$, $i=1,\cdots \mu$.
\end{prop}

\begin{preuve}
Pour toute $\omega =a\omega _{0}+b\omega _R$ de 
$\widehat\Omega ^{1}(\log S)$ on a, en tenant compte de la formule (\ref{derivation}),
\begin{eqnarray*}
\omega \in \widehat B^{1}_{DR}(\log S) &\Leftrightarrow &
\exists f,  \ g\in \widehat{\cal O}_2, \ 
\omega=df+g\cdot \omega _0 \\
&\Leftrightarrow &
\exists f,  \ g\in \widehat{\cal O}_2,\   
\omega=-X_0(f)\omega _{R}+ (R(f)+g)\omega _0\\
&\Leftrightarrow &
\exists f,  \ l\in \widehat{\cal O}_2,\   
\omega=X_0(f) \omega _{R}+ l\ \omega _0
\end{eqnarray*}
Ainsi l'obstruction à être un bord relatif ne porte que sur le seul coefficient de $\omega _R$
et se ramène à rechercher les obstructions à la résolution de l'équation
$X_0(f)=b$.~\end{preuve}

\bigskip

Cette proposition nous donne une formule intégrale pour le calcul des
coefficients de la forme normale de $X=X_0+\cdots$.
Pour cela, écrivons la version duale du théorème 2 énoncé dans l'introduction: 

\bigskip

{\bf Théorème $2^\sharp$.}
\textit{Soit $\omega=\omega _0 +\cdots $ une perturbation de la forme fermée  
$\omega _0$. Il existe un 
élément $(d_1,\cdots d_{\mu })$ de $\C[[h]]^{\mu}$, 
une conjugaison formelle fibrée $\Phi $ conjuguant orbitalement $\omega$ à
$$\omega_N=\omega _0 +\sum_{i=1}^{\mu }d_i(h){a_i\omega_R}.$$}

\bigskip

Les 1-formes apparaissant dans cette forme normale $\omega_N$ étant
relativement fermées par rapport à $\omega_0$,
il est naturel de les intégrer sur des cycles dans les feuilles de $\omega_0$,
c'est-à-dire dans les fibres de $h$. Pour ceci la base des formes
${a_i\omega_R}$
n'est pas adaptée: une base duale de cycles pour une fibre donnée ne l'est
plus pour les autres fibres lorsqu'on suit ces cycles par trivialisation locale
du fibré de Milnor induit par $h$. Nous choisissons donc une base de formes $\eta_i$
\textit{horizontales} c'est-à-dire de formes $b_i\omega_R$ dont les coefficients sont
solutions de $\nabla(b)=0.$ D'après $(\ref{solutions})$, il suffit de prendre
$$\eta_i=h^{-r_i}a_i\omega_R,\ \mbox{avec} \ r_i=\frac{p_1k_i+p_2l_i}{\delta }.$$
Fixons une fibre $F$ de $h$ et une base $\gamma_1,\cdots \gamma_\mu$ de
l'homologie de $F$ duale de la restriction des formes $\eta_i$ à $F$. 
Soit $\gamma_i(h)$ les cycles obtenus sur les autres fibres de $h$ à partir de
la fibre $F$ par les
trivialisations locales de la fibration $h$.
D'après la formule (voir \cite{MAL}):
$$\frac{d}{dh}\int_\gamma (h) \eta= \int_\gamma \nabla \eta$$
ces formes horizontales restreintes aux fibres voisines restent duales des
cycles $\gamma_i(h)$. Ecrivons maintenant $\omega_N$
dans cette base horizontale:
$$\omega_N=\omega _0 +\sum_{i=1}^{\mu }\delta_i(h){\eta_i}.$$
Nous avons obtenu:
$$\delta_i(h)=\int_{\gamma_i(h)}\omega_N.$$
Remarquons que les coefficients $\delta_i(h)$ sont multivalués, uniformisables
sur un revêtement fini -les exposants $r_i$ sont rationnels- du disque image
de $h$.
Observons cependant que, en tant que \textit{fonctions}, les $\delta_i$ ne sont pas des invariants de la classe
formelle du feuilletage défini par $\omega$. Plus précisément,
les fonctions
$$\widetilde\delta_i(h)=\int_{\gamma_i(h)}\omega$$
n'ont pas de raison d'être conjuguées aux
fonctions $\delta_i(h)$ définies sur la forme normale $\omega_N$.
La véritable nature
de ces invariants apparait dans le paragraphe qui suit.

\section{Réduction finale.}

D'après la remarque (\ref{unicite.3}), les conjugaisons formelles fibrées $\Phi $ de la
forme
$\exp b(h)R$ où $b$ est une série formelle d'une variable, agissent sur l'ensemble des formes
prénormales formelles du feuilletage défini par $X=X _0+\cdots$. Il s'agit ici de
définir un représentant privilégié unique dans chaque orbite de cette action. Explicitons
celle-ci.

\begin{lem}\label{c.fibree} Soit $R= p_1x\partial_x+p_2y\partial_y$ et $\Phi $ un
difféomorphisme formel, tangent à l'identité,  à l'origine de $\C^2$. Les propriétés
suivantes sont équivalentes:

i- $\exists $b$ \in \widehat{\cal O}_1,\ \Phi =\exp b(h)R,$

ii- $\exists $u$ \in \widehat{\cal O}_1, \ u(0)=1,\ \Phi =(xu(h)^{p_1},yu(h)^{p_2}),$

iii- $\exists \varphi  \in \widehat{\cal O}_1,\ \varphi (0)=0, \ \varphi '(0)=1,\ 
h\circ \Phi =\varphi \circ h.$
\end{lem}

\begin{preuve}
Pour vérifier l'équivalence entre (i-) et (ii-), on considère l'intégrale première
méromorphe $F=x^{p_2}y^{-p_1}$ du champ radial $R$. Le fait que $\Phi $ soit fibrée ($\Phi
=\exp B.R$) équivaut à $F\circ \Phi =F$. En écrivant $\Phi $ sous forme $(xV, yW)$ où $V$
et $W$ sont des unités de $\widehat{\cal O}_2$, cette condition donne
$V^{p_2}=W^{p_1}$ d'où l'existence d'une unité $U$ telle que $V=U^{p_1}$ et $W=U^{p_2}.$
De plus, le développement formel 
$$\exp B.R= \sum_{i=0}^{+\infty} \frac{(B.R) ^{(i)}}{i!}$$
prouve en l'appliquant à $x$ et $y$ que $B=b(h)$ si et seulement si $U=u(h).$
L'équivalence entre (ii-) et (iii-) se déduit de la relation de quasi-homogénéité
$$h(xu^{p_1},yu^{p_2})=u^{\delta }h(x,y)$$
en posant $\varphi (h)=hu(h)^{\delta }.$
\end{preuve}

\begin{lem} Soit $\Phi $ une conjugaison vérifiant une des conditions du lemme (\ref{c.fibree}).
On a:
\begin{eqnarray*}
\Phi ^*R&=&\frac{\varphi (h)}{h\varphi '(h)}\ R,\\
\Phi ^*X_0&=&(\frac{\varphi (h)}{h})^{\delta _0/\delta }\ X_0.
\end{eqnarray*}
\end{lem}

\begin{preuve}
Soit $\omega _0=\delta ^{-1}dh/h.$ On a
$$\Phi ^*\omega _0=\frac{h\varphi '(h)}{\varphi (h)}\ \omega _0.$$
De plus, le développement formel de $\Phi =\exp b(h)R$ montre l'existence d'un coefficient
$c(h)$ tel que $\Phi ^*R=c(h)R.$ De la relation $\omega _0(R)=1$, nous déduisons
$$c(h)=\frac{\varphi (h)}{h\varphi '(h)}.$$
Posons $\Phi ^*X_0=aX_0+bR.$ De la relation $\omega _0(X_0)=0$, nous déduisons $b=0$. 
Le développement formel de $\Phi =\exp b(h)R$ montre que $a=a(h)$. On calcule maintenant ce
coefficient en écrivant
$$[\Phi ^*R,\Phi ^*X_0]=\Phi ^*[R,X_0]=\Phi ^*\delta _0X_0=\delta _0a(h)X_0$$
soit encore
$$[\frac{\varphi (h)}{h\varphi '(h)}\cdot R, a(h)X_0]=\delta _0a(h)X_0$$
d'où
$$\frac{\varphi (h)}{h\varphi '(h)}(\delta a'(h)h+\delta _0a(h))=\delta _0a(h).$$
Cette équation différentielle d'une variable $h$ et d'inconnue $a$ admet pour solution 
$$a(h)= (\frac{\varphi (h)}{h})^{\delta _0/\delta }$$
d'où le résultat.
\end{preuve}

\bigskip

Considérons maintenant une forme prénormale $Y$ du feuilletage défini
par $X=X_0+\cdots$ sous la forme obtenue au théorème 1:
$$Y=X_0 +\sum_{i=1}^{\mu } d_i(h)a_iR.$$
On suppose qu'on a choisi une base $a_1=1,\ a_2,\cdots a_{\mu}$ de $\widehat{\cal O}_2 
/ J_0$ formée de monômes $a_i=x^{k_i}y^{l_i}.$ Ceci est toujours possible d'après
\cite{BG}. Pour toute conjugaison $\Phi $ vérifiant une des conditions du lemme
(\ref{c.fibree}), on a donc
$$\Phi ^*a_i=a_i(xu(h)^{p_1},yu(h)^{p_2})=u(h)^{k_ip_1+l_ip_2}a_i=
(\frac{\varphi (h)}{h})^{r _i}a_i$$
où $r _i$ est le rationnel $\frac{k_ip_1+l_ip_2}{\delta }$.
Nous obtenons donc:
\begin{eqnarray*}
\Phi ^*Y&=& (\frac{\varphi (h)}{h})^{\delta _0/\delta }\ X_0
+\sum_{i=1}^{\mu }d_i((\varphi (h))\frac{\varphi (h)^{1+r _i}}{h^{1+r _i}\varphi
'(h)}\ a_iR \\
&\approx & X_0 +\sum_{i=1}^{\mu } d_i((\varphi (h))
\frac{\varphi (h)^{1+q_i/\delta }}{h^{1+q_i/\delta }\varphi '(h)}\ a_iR\\
\end{eqnarray*}
où $q_i$ est l'entier $k_ip_1+l_ip_2-\delta _0$. Nous remarquons donc que si on
pose $\delta _i(h)=h^{1+q_i/\delta }d _i(h)$, le coefficient $\delta _i(h)$ change sous l'action
de $\Phi $ par la formule
$$\delta _i(h)\mapsto \frac{\delta _i\circ \varphi }{\varphi '}(h)$$
et se comporte donc comme un champ de vecteurs (ramifié) d'une variable $z$: 
$$\theta _i=z^{1+q_i/\delta }d_i(z)\partial _z.$$
Posons $d_i(z)=z^{m_i}u_i(z)$, avec $u_i(0)\neq 0$.
Ce champ ramifié se relève par $z=\tilde{z}^\delta $ en un champ uniforme
$$\widetilde{\theta _i}=\delta  \tilde{z}^{1+q_i+\delta
m_i}u_i(\tilde{z}^\delta )\partial_{\tilde{z}}.$$ Celui-ci se normalise (voir par exemple
\cite{MR2}) par une transformation $\widetilde{\varphi }$ de la même classe (de convergence,
divergence, sommabilité,...) que $d_i$ en
$$\frac{\delta \tilde{z}^{1+q_i+\delta m_i}}{1+\lambda  \tilde{z}^{q_i+\delta
m_i}}\partial_{\tilde{z}}$$ 
Cette forme normale redescend en
$$\frac{ {z}^{1+q_i/\delta +m_i}}{1+\lambda {z}^{q_i/\delta +m_i}}\partial_{z}.$$
Nous avons donc montré l'existence d'une conjuguante $\Phi  $ qui normalise 
$d_i(h)$
sous la forme
$$\frac{ {h}^{m_i}}{1+\lambda  {h}^{m_i+q_i/\delta }}.$$
Remarquons que le nombre complexe
$\lambda  $ est le résidu de la forme 
$$\frac{u_i^{-1}(\tilde{z}^\delta )}{\tilde{z}^{1+q_i+\delta m_i}}d\tilde{z}.$$
Pour qu'il soit non nul, il est nécessaire que $\delta $ divise $q_i$. En conséquence, nous
avons démontré l'existence d'une transformation $\Phi $ agissant sur l'ensemble des formes
prénormales de $X$ et normalisant un des coefficients $d_i$ choisi arbitrairement sous forme
$$\frac{h^{m_i}}{1+\lambda  h^{m_i+n_i}}$$
avec $m_i$ et $n_i$ entiers, $m_i\geq 1$, $n_i\geq 0$, et $\lambda $ nul dès que $\delta $ ne
divise pas
$q_i=k_ip_1+l_ip_2-\delta _0$. Nous avons obtenu le

\begin{thm}[Classification formelle des feuilletages]\label{f.normale}
Soit $X=X_0 +\cdots $ une perturbation à séparatrices fixées du champ quasi-\-ho\-mo\-gène
$X_0=(dh/\delta h)^{\flat }=\delta^{-1}X_h$, de degré $\delta_0$,
$\mu$ le nombre de Milnor
de $X_0$ (ou de $h$),
et $a_i=x^{k_i}y^{l_i}$ une base monomiale de $\widehat{\cal O}_2  / J_0$ où $J_0$ désigne
l'idéal engendré par les composantes de $X_0$. Il existe un 
élément $( d_1,\cdots
 d_{\mu })$ de
$\C[[h]]^{\mu}$ et une conjugaison fibrée $\Phi $ tangente à l'identité tels que

i- $\Phi $ conjugue pour la relation $\approx$ le champ $X$ à $$X_0 +\sum_{i=1}^{\mu
} d_i(h)a_i R$$

ii- un des coefficients formels $d_i$, choisi arbitrairement parmi les coefficients $d_j$ non
nuls --par exemple le premier non nul de la suite ordonnée par le choix de la base $a_i$--
s'écrive sous la forme rationnelle
$$\frac{h^{m_i}}{1+\lambda  h^{m_i+n_i}}$$
avec $m_i$ et $n_i$ entiers, $m_i\geq 1$, $n_i\geq 0$. De plus, le coefficient
$\lambda $ est nul dès que
le degré $\delta$ de $h$ ne
divise pas
$k_ip_1+l_ip_2-\delta _0$.
\end{thm}

Quelques commentaires sur cet énoncé:

\medskip
 
i-) La transformation utilisée dans cette réduction finale est dans la même classe
de sommabilité que la forme prénormale considérée. En particulier, elle converge dès que
celle-ci est convergente. Les problèmes de divergence des formes normales
obtenues ne se posent donc que sur l'étape de prénormalisation.

\medskip

ii-) De cette réduction finale des feuilletages définis par
$X=X_0+\cdots$,
nous pouvons remarquer que nous avons
caractérisé chaque classe formelle de feuilletage par une \textit{collection de
$\mu$ champs de vecteurs d'une seule variable, définis à conjugaison commune
près sur un revêtement fini d'ordre $\delta=\deg h$ du disque image de
l'intégrale première $h$ de $x_0$}. Une interprétation
géométrique de l'algèbre de Lie engendrée par ces invariants s'impose donc, ce
que nous nous proposons de développer dans un travail ultérieur.

\medskip

iii-) L'unicité des formes normales et conjuguantes données par le théorème 
\ref{f.normale} n'est pas complète: Nous pouvons encore faire agir le
groupe des transformations fibrées $\exp b(h)R$ qui préservent le coefficient
normalisé sous forme $h^{m_i}/(1+\lambda  h^{m_i+n_i})$, c'est-à-dire du
groupe des difféomorphismes en $h$ qui se relèvent en $\tilde h=h^\delta$ en un
difféomorphisme préservant le champ $\theta_i$. Ce groupe est une extension
abélienne du  groupe à un paramètre de $\theta_i$ par un groupe fini de
rotations. Son action modifiera les $\mu -1$ autres coefficients $d_j(h)$. Il n'est cependant pas nécessaire de choisir un représentant
privilégié sous cette action pour déterminer si deux champs $X=X_0+\cdots$
et $X'=X_0+\cdots$ sont formellement conjugués. La remarque ii) ci-dessus
répond à cette question en la réduisant à un problème de classification
d'algèbre de Lie d'une seule variable, et à ce titre, ce théorème peut donc
être considéré comme un théorème de classification complète.

\bigskip

Appliquons cette réduction finale des feuilletages aux trois exemples exposés dans les
paragraphes 2 et 3.

\bigskip

{\bf Exemple 1.} (Poincaré-Dulac) 
Toute forme prénormale formelle de
$\frac{1}{2pq}(qx\partial_x-py\partial_y)+\cdots$  est du type
$$X_0+d[[x^py^q]]R.$$
La normalisation finale donne ici
$$\Phi ^*X_0=X_0,\ \ \Phi ^*R=\frac{\varphi (h)}{h\varphi '(h)}\ R$$
et transforme $d(h)$ en $d\circ \varphi (h) \frac{\varphi (h)}{h\varphi '(h)}$.
Le coefficient $\delta (h)=hd(h)$ se normalise comme un champ de vecteurs, d'où la forme
normale
$$X_0+ \frac{h^m}{1+\lambda h^m} R,$$ encore équivalente à
$$(1+\lambda  h^m)\ X_0+ h^m\ R.$$
Celle-ci est une variante des formes normales formelles usuelles: pour retrouver les formes
normales $X_{p/q,m,\lambda }$ proposées dans \cite{MR2}, il suffit de faire une translation
de $1/2$ sur le résidu, et de faire agir un changement de variable linéaire convenable.

\bigskip

{\bf Exemple 2.} (singularité de type noeud-col).
Les formes prénormales formelles de $y\partial_y +\cdots $ sont ici du type
$$y\partial_y +{d}(x)x\partial _x.$$
Les transformations considérées dans la réduction
finale ($\Phi =\exp b(x)\partial_x$: $(x,y)\mapsto (\varphi(x),y)$) vérifient ici:
$$\Phi ^*X_0=X_0, \ \ \Phi ^*R=\frac{\varphi }{x\varphi '(x)}R$$
où $\varphi (x)=\exp b(x)\partial_x(x).$
En posant $d(x)R=\delta (x)\partial_x$,
une telle conjugaison agit sur $\delta (x)$ par $\frac{\delta \circ \varphi }{\varphi '}(x)$ et
permet de le normaliser sous la forme normale d'un champ. On
obtient donc la forme normale formelle usuelle
$$X\approx y\partial_y + \frac{x^{m+1}}{1+\mu  x^m}\partial_x
\approx {(1+\mu  x^m)}y\partial_y + {x^{m+1}}\partial_x.$$

{\bf Exemple 3.} (Le cusp $h=y^p-x^q$.) 
Les formes prénormales sont ici
$$X_0 +
%\sum_{\begin{array}{c}{\scriptstyle k=0,\cdots q-2,}\\ {\scriptstyle l=0\cdots p-2}\end{array}}
\sum_{k=0,\cdots q-2,\ l=0\cdots p-2}
 d_{k,l}(h)x^ky^l\ R.$$ D'après le théorème \ref{f.normale}, la transformation finale permet de normaliser un des coefficients $d_{k,l}(h)$ 
sous forme 
$$\frac{h^{m}}{1+\lambda  h^{m+n}}$$
le résidu $\lambda $ étant nul dès que
$\delta=p\vee q$ ne divise pas $(k_i+1)p_1+(l_i+1)q_1$.
On retrouve ainsi les formes normales formelles proposées par Frank Loray dans \cite{LOR}.

\bigskip

Montrons enfin comment nous pouvons déterminer les formes normales
formelles des \textit{champs} $X=X_0+\cdots$ lorsque le degré $\delta_0$
de $X_0$ est nul, à partir des formes prénormales obtenues au point 2) du
théorème 1. Examinons d'abord le cas du noeud-col (exemple 2). Les
formes prénormales formelles obtenues au paragraphe 3 sont:
$$X\sim X_1=a_1(x)X_0+b_1(x)R, \ X_0=y\partial_y,\ R=x\partial_x,\ \nu(b_1)=m>0.$$
D'après la remarque $\ref{unicite.1}$, on peut faire opérer sur ces formes
prénormales de $X$ tout difféomorphisme formel $\Phi =\exp Z$
tel que $Z$ commute avec $X_0$, c'est-à-dire, puisque $[X_0,R]=0$, de la forme
$Z=\alpha (x)X_0+\beta (x)R.$
Considérons d'abord l'action d'un difféomorphisme de type $\exp\beta (x)R
$
c'est-à-dire d'un difféomorphisme $\Phi$ ne dépendant que de la seule variable
$x$:  $\Phi(x,y)=\varphi(x)$. La seconde composante de $X_1$ ne dépendant que
de cette seule variable, nous pouvons choisir $\Phi$ de sorte que 
$$X\sim X_2=a_2(x)X_0+b_2(x)R,\ \mbox{avec }b_2(x)=\frac{x^m}{1+\lambda x^m},
\ a_2=a_1\circ\varphi .$$
On peut maintenant faire agir une conjugaison formelle de la forme $\Phi=
\exp\alpha (x)X_0.$ Elle laisse invariante la première composante de $X_2$
puisque $\alpha (x)X_0$ commute avec $a_2(x)X_0$, et agit sur la seconde par:
\begin{eqnarray*}
(\exp\alpha (x)X_0)^*\ b(x)R&=&
b(x)R + [\alpha (x)X_0,b(x)R ]+ \\
&+&\frac{1}{2}[\alpha (x)X_0,[\alpha (x)X_0,b(x)R ]]+\cdots\\
&=&
b(x)R -b(x)x\alpha'(x)X_0    .
\end{eqnarray*}
d'où
$$(\exp\alpha (x)X_0)^*\ X_2
= (a_2(x)-b_2(x)x\alpha'(x))X_0 + b_2(x)R.$$
Posons: $a_2(x)=P_m(x)+x^{m+1}v(x)$ où $P_m$ est un polynôme de degré au plus
$m$, et choisissons $\alpha(x)$ de sorte que 
$$v(x)-\frac{\alpha '(x)}{1+\lambda x^m}=0.$$
Nous obtenons
$$X\sim P_m(x)y\partial_y +\frac{x^m}{1+\lambda x^m}x\partial_x.$$
Ces formes normales sont celles obtenus par A.D. Bruno \cite{BRUNO}, et sous
une variante par L. Teyssier \cite{TEY1}. Ce dernier a poursuivi la
classification analytique de ces champs dans \cite{TEY2}.

\bigskip

Un raisonnement identique (la première étape étant déjà explicitée lors de la
réduction finale du feuilletage) donne dans le cas des champs de Poincaré-Dulac
(exemple 1) les formes normales formelles suivantes:
$$X\sim P_m(h)X_0 +\frac{h^m}{1+\lambda h^m}R,$$
où $P_m$ est un polynôme de degré au plus
$m$.

\bigskip
\noindent {\sc E. Paul.}
\medskip

\noindent Laboratoire Emile Picard, U.M.R. C.N.R.S. 5580.

\noindent Université Paul Sabatier, 118 route de Narbonne, 31062 Toulouse Cedex 4, France.

\noindent email: paul@picard.ups-tlse.fr


\begin{thebibliography}{99}

\bibitem{BG} {\sc J. Briançon, A. Galligo.} 
{\em Déformations distinguées d'un point de $\C^2$ ou $\R^2$,} dans
Singularités à Cargèse, Astérisque 7 et 8, (1973) p 129--138.

\bibitem{BRUNO} {\sc A. D. Bruno             } 
{\em Local methods in nonlinear differential equations     ,} 
Springer-Verlag (1989)
%\bibitem{CMA} {\sc D. Cerveau, J.F. Mattei,}
%{\em Formes intégrables holomorphes singulières.}
%Astérisque S.M.F.,
%{\bf 97}, (1982).

\bibitem{C-D}{\sc M. Canalis-Durand, F. Michel, M. Teisseyre.}
{\em  Algorithms for formal reduction of vector field singularities.}
 J. Dynam. Control Systems {\bf 7, 1}  (2001), p. 101--125. 

\bibitem{CM} {\sc D. Cerveau, R. Moussu.}
{\em Groupes d'automorphismes de $(\C,0)$ et équations diffé\-ren\-tielles $y
\ dy+\cdots =0.$ }
Bull. Soc. math. France, {\bf 116}, (1988), p.459--488.

\bibitem{LOR} {\sc F. Loray.} 
{\em Réduction formelle des singularités cuspidales de champs de vecteurs
analytiques.} J. of Diff. Equations {\bf 158, 1} (1999), p. 152--173.

\bibitem{MAL} {\sc B. Malgrange.} 
{\em Intégrales asymptotiques et monodromie.}
Ann. scient. Ec. Norm. Sup.  série 4, t.{\bf 7}, (1974) p. 405--430.

\bibitem{MAL.2} {\sc B. Malgrange.} 
{\em Sur les points singuliers des équations différentielles.}
L'enseignement mathématique, T. XX, fasc. 1-2, (1974) p. 147--176.

\bibitem{MAR} {\sc J. Martinet.} 
{\em Normalisation des champs de vecteurs holomorphes.}
Séminaire Bourbaki 1980-1981, exposé 564, Lecture notes in Math. {\bf 55-70} (1981).

\bibitem{MR2} {\sc J. Martinet, J.P. Ramis.}
{\em Classification analytique des équations diffé\-ren\-tielles non
linéaires résonnantes du pemier ordre.}
Ann. Sci. Ecole Norm. Sup.,
t. {\bf 16}, (1983), p. 571--621.

\bibitem{JFM} {\sc J.F. Mattei.}
{\em Quasihomogénéité et équiréductibilité de feuilletages holomorphes en dimension 2.}
Astérisque {\bf 261} (2000) p. 253--276.

\bibitem{MM} {\sc J.F. Mattei, R. Moussu.}
{\em Holonomie et intégrales premières.}
Ann. Sci. Ecole Norm. Sup.,
t.{\bf 13}, (1980), p. 469--523.

\bibitem{SAITO.1} {\sc K. Saito.}
{\em Quasi homogene isolierte singularitaten von hyperflachen,}
Inven. Math. {\bf 14}, p.123--142 (1971).

\bibitem{SAITO.2} {\sc K. Saito.}
{\em Theory of logarithmic differential forms and logarithmic vector fields.}
Journal of the Faculty of Sciences of Tokyo,
vol. {\bf 27, 2} (1980). 

%\bibitem{SCD} {\sc R. Schaefke, M. Canalis-Durand.} ???

\bibitem{STOLO} {\sc L. Stolovitch.}
{\em Singular complete integrability.}
Pub. Math. I.H.E.S. {\bf 91}, (2000).

\bibitem{TEY1} {\sc L. Teyssier.}
{\em Equations homologique et cycles asymptotiques d'une singularité
  noeud-col,}
Preprint IRMA Lille, vol {\bf 55}, (2001)

\bibitem{TEY2} {\sc L. Teyssier.}
{\em Analytic classification of singular saddle-node vector fields,}
Preprint Rennes, 03-02 (2003).


\bibitem{ZOL} {\sc E. Str\'ozyna, H. Zoladek.} 
{\em The analytic normal form for the nilpotent singularity,}
J. of Diff. Equations {\bf 179, 2} (2002), p. 479--537.




\end{thebibliography}
\end{document}